\documentclass[11pt,psamsfonts]{amsart}
\usepackage{amsmath}
\usepackage{amsthm}
\usepackage{amssymb}
\usepackage{amscd}
\usepackage{amsfonts}
\usepackage{amsbsy}
\usepackage{epsfig,afterpage}
\usepackage[dvips]{psfrag}
\usepackage{color}     

\newcommand{\R}{\ensuremath{\mathbb{R}}}

\newcommand{\N}{\ensuremath{\mathbb{N}}}

\newcommand{\fp}{\hfill\blacksquare}
\newcommand{\dem}{\textbf{Proof}}
\newcommand{\dis}{\displaystyle}
\newcommand{\vs}{\vspace{0,5cm}}
\def\omegar{\Omega}

\newtheorem {theorem} {Theorem} 
\newtheorem {proposition} [theorem] {Proposition}

\newtheorem {lemma} [theorem] {Lemma}
\newtheorem {definition} [theorem] {Definition}
\newtheorem {remark} {Remark}

\begin{document}

\title[Two-fold symmetric singularity] {Two-fold symmetric singularity}

\author[A. Jacquemard, M.A. Teixeira, D.J. Tonon]
{Alain Jacquemard$^1$, Marco A. Teixeira$^2$, Durval J. Tonon$^3$}

\address{$^1$ Institut de Math\'{e}matiques de Bourgogne
UMR CNRS 5584, Universit\'{e} de Bourgogne, Dijon, France\\
$^2$ Department of Mathematics IMECC, Universidade Estadual de
Campinas Rua Sergio Buarque de Holanda 651, Cidade Universit´aria -
Bar˜ao Geraldo, 6065 Campinas (SP), Brazil\\
$^3$ Universidade Federal de Goi\'{a}s, IME, CEP 74001-970 – Caixa
Postal 131, Goi\^{a}nia, GO, Brazil.}

\email{jacmar@u-bourgogne.fr} \email{teixeira@ime.unicamp.br}
\email{djtonon@mat.ufg.br}


\date{}
\dedicatory{}

\maketitle

\begin{abstract}
We explore some qualitative dynamics in the neighborhood of the
$3-dimensional$ two-fold symmetric singularity. We study the
existence of an one-parameter family of regular (pseudo) periodic
orbits of such systems near a reversible two-fold singularity.
\end{abstract}


\section{Introduction}

Theory of non-smooth dynamical systems has been developing at a very
fast pace in recent years and it has become certainly one of the
common frontiers between Mathematics and Physics or Engineering.
Hybrid and switched models are being increasingly used in
applications to describe a large variety of physical devices.
Examples include mechanical systems with friction and backlash,
electrical and electronic circuits, walking and hopping robots and,
more recently, biological and neural systems \cite{B-S-C,C,K-R-G}.
One should observe that much work has been done in the study of the
qualitative aspects of the phase space of discontinuous vector
fields.

We consider vector fields expressible in the form $\dot x = Z(x)$
where $x \in \R^3$ is a state vector and $Z$ is a smooth piecewise
mapping. The discontinuities are concentrated on a codimension-one
submanifold $\Sigma$ of $\R^3$. $\Sigma$ is usually called the
switching manifold. Orbit-solutions on $\Sigma$, whenever possible
are defined according to the Filippov convention, \cite{F}. We point
out that trajectories may become constrained to the switching
manifold, and this behavior is called sliding.

A general understanding of dynamics of generic $3$D Filippov sytems
was obstructed by the appearance of the {\it two-fold singularity}
(refer to \ref{two-fold section} and \cite{T1}). Observe that,
topologically speaking, its shape  is very simple, and moreover it
is generic in piecewise smooth systems with three or more
dimensions.

In the last years we have been noticing a great interest in the
study of the two-fold singularity (also called T-singularity). In
\cite{J-C} is classified the dynamics in a neighborhood of a
T-singularity. In \cite{C-B-F-J} her occurrence was discussed with
real models in control theory, introducing conditions for the study
of the existence of this singularity.

Let us consider a discontinuous vector field $Z=(X,Y)$ and its
corresponding {\it normalized sliding vector field}
$\widetilde{Z}^S$ (see Definition \ref{sliding} and Remark
\ref{normalized-sliding}). Let $p_0$ a point where $Z$ presents a
discontinuity. When $p_0$ is the two-fold singularity, $p_0$ is
always a critical point of $\widetilde{Z}^S$ and a fixed point of
the {\it return map} (see definition
\ref{aplicacao-primeiro-retorno}) which is a composition
$\varphi_Z^{}=\gamma_Y^{} \circ \gamma_X^{}$ of two involutions
associated respectively to $Y$ and $X$.

We deal with perturbations $Z$ (whose expressions will be precise in
$\ref{normal-forms-section}$) of the $3D$ piecewise smooth vector
fields $Z_0=(X_0,Y_0)$. Our aim is to find conditions to the
existence of typical closed orbits of $Z$ (see definition
$\ref{definicao-orbitapseudoperiodica}$). To do this we use some
geometrical properties, like a reversibility. So, here we exhibited
a subset of this typical singularities that is reversible. 

The paper is organized as follows: in sections \ref{preliminares} we
deal with some preliminaries, give some definitions, and establish
the notations. In section \ref{reversibilidade} we state the results
about the reversibility. In section \ref{reversibilidade} and
\ref{pseudo-periodica} we prove the existence of a family of the
regular and pseudo periodic orbits, respectively.



\section{Preliminaries}\label{preliminares}

\subsection{Distinguished regions of the discontinuity set} 
\label{subsection:regions} In this section some notations, basic
definitions and elementary concepts are presented.

Designate by  $\mathfrak{X}^r$ the space of all germs of $C^r$
vector fields on ${\R}^{3}$ at $0$ endowed with the $C^r$--topology
with $r> 1$ and large enough for our purposes. Call
\textbf{$\Omega^r$} the space of all germs of vector fields $Z$ in
$(\R^{3},0)$ such that

\[
\label{discontinuity} Z(q)=\left\{\begin{array}{l} X(q),\quad $for$
\quad h(q)>0,\\ Y(q),\quad $for$ \quad h(q)<0,
\end{array}\right.
\]where $h$ is (a germ of) a smooth function $h:(\R^{3},0)
\rightarrow (\R,0)$ having $0 \in \R$ as its regular value. Let
$\Sigma=h^{-1}(0)$. We assume that $0 \in \Sigma.$

The above vector field is denoted by $Z = (X,Y)$. We endow
$\Omega^r=\mathfrak{X}^r\times \mathfrak{X}^r$  with the product
topology.

To define the orbit solutions of $Z$ on the switching surface
$\Sigma$ we follow a pragmatic approach. In a well characterized
open set $\mathcal{O}$ of $\Sigma$ (described below) the solution of
$Z$ through a point $p \in \mathcal{O}$ obeys the Filippov rules
(see \cite{F}) and on $\Sigma-\mathcal{O}$ we accept it to be
multivalued. As we are dealing with systems derived from ordinary
differential equations the non-uniqueness of solutions is allowed.
We just must take into account all the leaves of the foliation in
$\R^3$ generated by the orbits of $Z$ (and also the orbits of $X$
and $Y$) passing through or exiting from or converging to a point
$p\in\Sigma$.

For each $X \in \mathfrak{X}^r$ we define the smooth function
$Xh:\R^3\rightarrow\R$ given by $Xh=X.\nabla h$ where $.$ is the
canonical scalar product in $\R^3$.

\noindent In what follows we use the Filippov convention (refer to
\cite{F}). We first distinguish the following regions on $\Sigma$ :
\begin{enumerate}

\item [$\blacktriangleright$] \textbf{Sewing Region}: $SwR=\{p\in \Sigma;(Xh)(p)(Yh)(p)>0\}$. When
convenient we denote $SwR^+=\{p\in \Sigma;(Xh)(p)>0,(Yh)(p)>0\}$ and
$SwR^-=\{p\in \Sigma;(Xh)(p)<0,(Yh)(p)<0\}$. In general a point in
the phase space which moves on an orbit of $Z$ and reaches a point
in $SwR$, crosses $\Sigma$.

\item [$\blacktriangleright$] \textbf{Escaping Region}: $EscR=\{p\in \Sigma;(Xh)(p)>0 ,(Yh)(p)<0\}$. In this case
any orbit which meets $EscR$ remains tangent to $\Sigma$ for
negative times.

\item [$\blacktriangleright$] \textbf{Sliding Region}: $SlR=\{p\in \Sigma;(Xh)(p)<0,(Yh)(p)>0\}$.
On $SlR$ the flow slides on $\Sigma$; the flow follows a well
defined vector field $Z^S$ called the {\it sliding vector field}
(see Definition~\ref{sliding}).
\end{enumerate}

\noindent Generically, the set $\mathcal{O}= SlR \cup EscR \cup SwR$
is open and dense in $\Sigma$. Observe that for any $p \in
\mathcal{O}$, we have $X(p)\neq 0$ and $Y(p) \neq 0$.

\begin{definition} Let $Z=(X,Y) \in \Omega^r$. The sliding vector field $Z^S$ associated to $Z$ is a linear convex combination of $X$ and $Y$ tangent to
$\Sigma$, that is,
\[
Z^S=\frac{1}{(Y-X).\nabla h}(Y. \nabla hX - X.\nabla hY).
\]
\label{sliding}\end{definition}

\begin{remark} Observe that $EscR$ for $Z$ represent $SlR$ for $-Z$. We can define the {\it
escaping vector field} on $\Sigma$ by $-(-Z)^S$. This vector field
is called the \emph{sliding vector field} independently of whether
it is defined in the sliding or escaping region. For $p\in SlR\cup
EscR$ the local orbit of $p$ is ruled by this vector field.
Therefore all this orbit is contained in $SlR\cup EscR$, and the
future orbit of the sliding vector field coincides with the future
orbit of the normalized sliding vector field
$\widetilde{Z}^S=(Y.\nabla hX - X.\nabla hY)$.
\label{normalized-sliding}\end{remark}


\noindent Observe that $Z^S$ and $\widetilde{Z}^S$ are orbitally
equivalent on $SlR$ ({\it resp.} on $EscR$). Moreover
$\widetilde{Z}^S$ can be $C^r$-extended beyond the boundary of
$SlR,EscR$. For technical reasons we consider the future orbit of
$Z$ through a point $p\in \overline{SlR\cup EscR}$ given by the
orbit of $\widetilde{Z}^S$.

\vspace{0,5cm}

\noindent \textbf{Notation:} In all what follows we consider the map
$h:(x,y,z)\mapsto z$. So, the expression of $\widetilde{Z}^S$ is:

\[
\widetilde{Z}^S=(X^1Y^3 -Y^1X^3, X^2Y^3- Y^2X^3),
\]where
$X=(X^1,X^2,X^3)$ and $Y=(Y^1,Y^2,Y^3)$.


%

\begin{definition}
We say that $0$ is a {\it two-fold singularity} of $Z=(X,Y) \in
\Omega^r$ if $Xh(0)=Yh(0)=0$ and $X^2h(0)\neq 0, Y^2h(0)\neq0$,
where $X^2h(0)=X(Xh)(0)$.
\end{definition}

\begin{definition}
If $p \in SlR\cup EscR$ and $X(p),Y(p)$ are linearly dependent then
$p$ is a critical point of $Z^S$. In this case $p$ is called a
pseudo equilibrium of $Z$.
\end{definition}

\noindent The curves of {\it tangential singularities} or the
\emph{$\Sigma-$singularity} of $X$ in $\Sigma$ are given by
$S_X=\{p\in \Sigma;Xh(p)=0\}$.

\noindent We denote  by $\Omega^F\subset \Omega^r$ the set of
non-smooth vector fields such that the origin is a two-fold
singularity. First of all, observe that:

\begin{itemize}
\item[1-] The trajectories of both $X$ and $Y$ through $0$ have a quadratic
contact with $\Sigma$ (at $0$).

\item[2-] Generically, $S_X$ and $S_Y$ are transverse at $0$. This case
was first studied in \cite{T1}.

\item [3-] If $0$ is a two-fold singularity of $Z$ then $\widetilde{Z}^S$ can
be $C^r$-extended to a full neighborhood of
$0\in\Sigma,\widetilde{Z}^S(0)=0$ and $\varphi_Z(0)=0$.
\end{itemize}

\subsection{Two-fold singularity}\label{two-fold section}

The following construction is given in \cite{T2}. Let $Z=(X,Y) \in
\Omega^r$ such that $0$ is a fold point for {\it both} $X$ and $Y$.
Applying the Implicit Function Theorem, for each $p \in (\Sigma,0)$
there exists a unique $t(p)$ such that the orbit-solution $t\mapsto
\phi_X(t,p)$ of $X$ through $p$ meets $\Sigma$ at a point
$\widetilde{p}=\phi_X(t(p),p)$.

We then define the smooth mapping $\gamma_X:(\R^2,0)\rightarrow
(\R^2,0)$ by $\gamma_X(p)=\widetilde{p}$. This map is a
$C^r$-diffeomorphism and satisfies: $\gamma_X^2=Id$.
Analogously, we define the smooth map associated to $Y$:
$\gamma_Y:(\R^2,0)\rightarrow (\R^2,0)$ which satisfies
$\gamma_Y^2=Id$. We define now the first return map associated to
$Z=(X,Y)$:

\begin{definition}
The first return map $\varphi_Z: (\Sigma,0)\rightarrow (\Sigma,0)$
is defined by the composition
\[
\varphi_Z=\gamma_Y\circ \gamma_X
\]
We denote by $L_Z(.,.)$ the linear part of $\varphi_Z$.
\label{aplicacao-primeiro-retorno}\end{definition}

\noindent $\varphi_Z$ is a $C^r$-diffeomorphism preserving area:
$det(D\varphi_Z)(0,0)=1$. So the eigenvalues of $D\varphi(0,0)$ are
$\beta$ and $\beta^{-1}$ with:

\begin{itemize}
    \item [$(a)$] saddle type: $\beta\in \R$, $\beta\neq 0$;
    \item [$(b)$] elliptical type: $\beta=e^{i\,\theta}$ with $\theta\in ]0,\pi[$.  In this case, $L_Z(.,.)$ is a rotation.
\end{itemize}

%
%
%

\subsection{Partition of $\Omega^F$} Consider the subset of $\Omega^F$:

\begin{itemize}
    \item []{\bf Regular two-fold:} Let $\Omega_0^F(\delta)$ be the set of
all $Z=(X,Y)\in \Omega^F$ such that the contact between $S_X$ and
$S_Y$ at $0$ is transverse, the eigenvectors of
$D\widetilde{Z}^S(0)$ are transverse to $S_X$ and $S_Y$ at $0$ and
$0$ is a hyperbolic critical point for $\widetilde{Z}^S$.
\end{itemize}

\begin{remark}
From \cite{S-T1} we derive that $\Omega_0^F(\delta)$ is a
codimension zero subma\-ni\-fold of $\Omega$.
\end{remark}

\begin{figure}[ht]
\epsfysize=5cm
\psfrag{Sx}{$S_X$}\psfrag{Sy}{$S_Y$}\psfrag{M}{$\Sigma$}\psfrag{1}{Parabolic
case}\psfrag{2}{Elliptic case}\psfrag{3}{Hyperbolic case}
\psfrag{E}{$h_2>1$}\centerline{\epsfbox{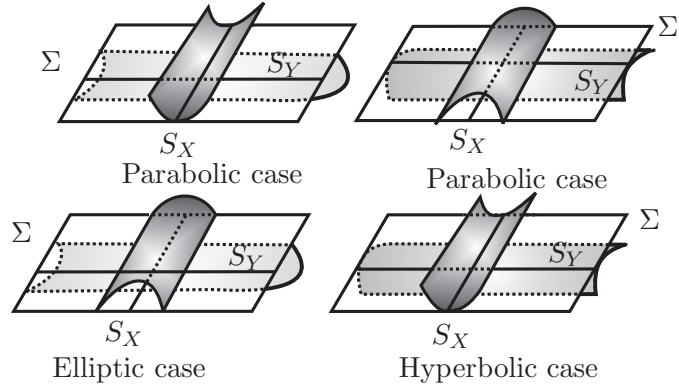}}\caption{Two-fold
singularity.}\label{Tiposdedobra}
\end{figure}

In $\Omega_0^F(\delta)$ we distinguish the following subsets,
corresponding to three cases:

\textbf{Elliptic case}: $\Omega_0^F(\delta.1)=\{Z\in
\Omega_0^F(\delta);X^2h(0)<0$ and $Y^2h(0)>0\}$. We have two
invisible tangencies (\emph{invisible two-fold} or \emph{a
T-singularity}). See Figure $\ref{Tiposdedobra}$;

\textbf{Parabolic case}: $\Omega_0^F(\delta.2)=\{Z\in
\Omega_0^F(\delta);X^2h(0)>0, Y^2h(0)>0$ or $X^2h(0)<0,Y^2h(0)<0\}$
(visible fold- invisible fold);

\textbf{Hyperbolic case}: $\Omega_0^F(\delta.3)=\{Z\in
\Omega_0^F(\delta);X^2h(0)>0, Y^2h(0)<0\}$ (visible two-fold).

Recently, many works tried to understand the dynamics of the
T-singularity, \cite{J-C, J, C-B-F-J}. See \cite{T2} for further
references and related topics.

\begin{definition} We say that a $C^{\infty}$-diffeomorphism
$\xi:\R^3\rightarrow \R^3$, is an involution if $\xi\circ \xi = Id$.
We put $Fix(\xi)=\{p \in \R^4\ | \ \xi(p)=p\}$.
\end{definition}

\begin{definition}
We say that a non-smooth vector field $Z\in \Omega^r$ with
discontinuity manifold $\Sigma$ is $\xi$-reversible if there exists
an involution $\xi:\R^3 \rightarrow \R^3$ such that:

\begin{enumerate}
\item[1-]  $Fix(\xi)\subset \Sigma$
\item[2-]  $\forall p \notin \Sigma$, $\xi\circ Z(p)=-Z(\xi(p))$.
\end{enumerate}
\label{aplicacao reversivel}
\end{definition}

\begin{definition}
Let $\Omega_0^{\xi}(\delta)$ be the set of elements
$Z\in\Omega_0^F(\delta)$ such that:
\begin{itemize}
    \item [1-] $Z$ is $\xi$-reversible;
    \item [2-] $Z$ is a generic invisible two-fold at 0;
\end{itemize}
\end{definition}

\begin{definition}
We say that $Z\in \Omega^r$ is simple at 0 if $X^2h(0)\neq \pm
XYh(0)$.
\end{definition}

We say that $X\in \mathfrak{X}^r$ is \emph{semi-linear} if
$X^1(x,y,z),X^2(x,y,z)$ possess only the 0-jet (degree zero) and
$X^3(x,y,z)$ possess only the 1-jet (degree one). Put:

\[
SL(\Omega_0^F(\delta))=\{Z=(X,Y)\in \Omega_0^F(\delta);X \mbox{ and
} Y \mbox{ are semi-linear}\}.
\]

%

Consider the set of non smooth dynamics vector fields:

\[
\begin{array}{ll} \Omega_0^P(\delta)                 &=\{Z\in SL(\Omega_0^F(\delta));X^2h(0)<0,(X(Yh))(0)<0,\\\\
                                                     &(X(Yh))(0)(Y(Xh))(0)=(X^2h)(0)(Y^2h)(0),Y^2h(0)>0\}.
\end{array}
\]


\subsection{Orbits}\label{orbits}

By convention, if $p \in SlR$ the future orbit of $Z$ through $p$ is
given by the trajectory of the sliding vector field $Z^S$ through
$p$.

%

\begin{figure}[ht]
\epsfysize=3cm
\psfrag{X}{$X$}\psfrag{1}{$(a)$}\psfrag{2}{$(b)$}\psfrag{Y}{$Y$}\psfrag{Sx}{$S_X$}\psfrag{Sy}{$S_Y$}
\psfrag{Rcu}{$SwR^+$}\psfrag{Rcd}{$SwR^-$}
\psfrag{F1}{$\varphi_1$}\psfrag{Rd}{$SlR$}\psfrag{F2}{$\varphi_2$}\psfrag{P}{$p_0$}\psfrag{Re}{$EscR$}
\psfrag{A}{$\phi_X$}\psfrag{B}{$\phi_Y$}\psfrag{M}{$\Sigma$}\psfrag{1}{Regular}\psfrag{2}{Pseudo}\psfrag{Q}{$q$}
\psfrag{P}{$p$}\psfrag{C}{$\phi_{Z^S}$}\centerline{\epsfbox{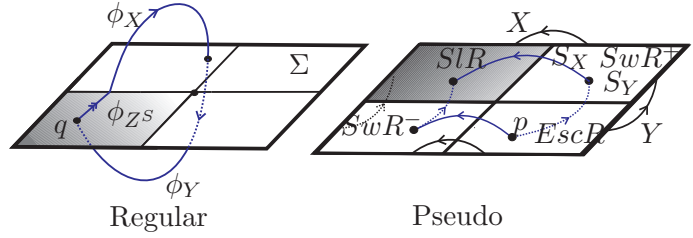}}
\caption{Periodic orbits.}\label{resultados}
\end{figure}

\begin{definition}Let $Z=(X,Y)\in \Omega^r$ and $x_0\in (\R^3,0)$.
\begin{itemize}
    \item [1-] The orbit $t\mapsto \phi_Z(t,x_0)$ is a regular periodic orbit if it is
closed and composed by segments of orbits of $X, Y$ and $Z^S$
keeping the orientation. See Figure $\ref{resultados}$.

    \item [2-] We call pseudo periodic regular orbit a trajectory of $Z$ which is
closed, composed by segments of the orbits of $X$ and $Y$, with non
preserved orientation (see Figure $\ref{resultados}$).
\end{itemize}
\label{definicao-orbitapseudoperiodica}\end{definition}

\subsection{Main results} Consider the set of straights line passing for
the origin:

\[
\mathcal{C}_{\alpha}=\{(x,\alpha x,0)\in \Sigma; \mbox{ with }\alpha
\in \R\}.
\]

We obtain the result about the existence of pseudo periodic orbits
and periodic orbits for $Z\in \Omega_0^{\xi}(\delta)$:

\begin{theorem}
Let $Z\in \Omega_0^{\xi}(\delta)$ and $p\in \mathcal{C}_{\alpha}$.
Then:

\begin{itemize}
    \item [$(a)$] If $Z$ is simple at 0 then $Z$ does not have k-periodic orbits for all $k\in \N$.

    \item [$(b)$] If $Z$ is non simple at 0 and $\alpha=1$ then $Z$ has one-periodic
    orbit passing through $p$. In addition, for $\alpha\neq 1$ does not exist $k$-periodic orbits for all
    $k\in \N$.

    \item [$(c)$] If $Z$ is non simple at 0 and $\alpha=-1$ then $Z$ has one-pseudo periodic
    orbit passing through $p$. In addition, for $\alpha\neq -1$ does not
    exist $k$-pseudo periodic orbits for all $k\in \N$.

     \item [$(d)$] 0 is a hyperbolic equilibrium saddle (elliptical) point of
    $\varphi_Z$ provided that $|XYh(0)|>|X^2h(0)|, (|XYh(0)|<|X^2h(0)|)$, see Figure
    \ref{tipo-primeiro-retorno}.
\end{itemize}\label{existenciadeciclos}\end{theorem}

In Theorem $\ref{existenciadeciclos}$ we exhibit conditions on $Z\in
\Omega_0^{\xi}(\delta)$ for the existence of one-parameter families
of periodic and pseudo periodic orbits for $Z$.

%
%
%
%
%
%
%
%
%
%
%
%
%

We also study the existence of pseudo periodic orbits and regular
periodic orbits for $Z\in SL(\Omega_0^F(\delta))$. Observe that
hypothesis in Theorem $\ref{proposicao-orbita-pseudo-periodica}$
$(a)$ complete the result obtained in Theorem
$\ref{existenciadeciclos}$ $(c)$ for $Z\in SL(\Omega_0^F(\delta))$.
In the item $(b)$ we get a more general family of regular periodic
orbits, generalizing the family obtained in Theorem
$\ref{existenciadeciclos}$ $(b)$,  for $Z\in
SL(\Omega_0^F(\delta))$. There results are summarized in:

\begin{theorem}
Consider $Z=(X,Y)\in SL(\Omega_0^F(\delta)),p\in
\mathcal{C}_{\alpha}$ and $X^2h(0)<0,Y^2h(0)>0$. Then
\begin{itemize}
    \item [$(a)$] If $Z$ is simple at the origin,
    $XYh(0)=-YXh(0),X^2h(0)=-Y^2h(0)$ and $\alpha=-1$ then $Z$ has a 1-pseudo periodic
    orbit passing through $p$. In addition, for $\alpha\neq -1$ there does not exist $k$-pseudo periodic orbits for all $k\in \N$.

    \item [$(b)$] If $Z\in\Omega_0^P(\delta), \alpha=\lambda^*$ and $p\in SwR$ then $p$ is either a fixed point for
    $\varphi_1$ or $\varphi_2$, where
$\lambda^*=\frac{X^2h(0)}{XYh(0)}$.
\end{itemize}\label{proposicao-orbita-pseudo-periodica}\end{theorem}


\section{Reversibility}\label{reversibilidade}

In this section we will explore some symmetry properties of $Z\in
\Omega_0^F(\delta)$. Our objective is to study the existence of
periodic orbits of arbitrary period $k$.

\subsection{Normal forms for the reversible two-fold singularity}\label{normal-forms-section} Consider $Z_0\in
\Omega_0^F(\delta)$. Throughout this section, we fix a coordinates
system $(x,y,z)$ such that

\[
\xi(x,y,z)=(y,x,-z). \] Observe that $Fix(\xi)=\{(x,y,z);
x=y,z=0\}$.

Let $f:\R^3\rightarrow\R$ be a non-zero polynomial in $(x,y,z)$. We
call $d_T{}(f)$ the {\it total degree} of $f$, that is the maximum
sum of exponents of the monomials of $f$.

Consider the polynomials: $f_{\sigma}^i:\R^3\rightarrow \R$ where
$\sigma\in\{x,y\}$ and $i=1,2,3$, with the following specifications:
$f_{\sigma}^{1,2}$ have null linear part and
$d_T{}(f_{\sigma}^3)>1$. Let

\[
\begin{array}{ll} F_{\sigma}(x,y,z)
&=(f_{\sigma}^1(x,y,z),f_{\sigma}^2(x,y,z),f_{\sigma}^3(x,y,z)).
\end{array}
\]

\noindent In what follows we exhibit the topological normal form for
the two-fold singularity:

\begin{proposition} (\textbf{Normal Forms}) If $Z \in \Omega_0^F(\delta)$ then $Z$ is $C^0-$equivalent to
$\widetilde{Z}=(\widetilde{X},\widetilde{Y})$ where

\[
\begin{array}{ll} \widetilde{X}(x,y,z)    &=(C_X,C_{XY},x)+F_x(x,y,z),\\
                  \widetilde{Y}(x,y,z)    &=(C_{YX},C_Y,y)+F_y(x,y,z),
\end{array}
\] with $C_X=X^2h(0),C_Y=Y^2h(0),
C_{XY}=X(Yh)(0),$ $C_{YX}=Y(Xh)(0),$ $X^2h(0)\neq 0$ and
$Y^2h(0)\neq 0$. \label{normal-forms}\end{proposition}

\dem.  The origin belongs to the frontier of the sliding region.
Consider a local coordinates system such that
$S_X=\{(x,y,z);x=z=0\}$ and $S_Y=\{(x,y,z);y=z=0\}$. With these
settings the local normal forms of $X$ and $Y$ are:

\[
\begin{array}{ll}   X(x,y,z)      &=(C_X, C_{XY},x)+F_x(x,y,z)\\\\
                   Y(x,y,z)       &=(C_{YX},C_Y,y)+F_y(x,y,z).
\end{array}
\]
Observe that $X^2h(0)=C_X,Y^2h(0)=C_Y,C_{XY}=X(Yh)(0)$ and
$C_{YX}=Y(Xh)(0)$.

$\fp$

\noindent Consider the regions in $\Sigma$:

\begin{equation}
\begin{array}{ll} SlR           &=\{(x,y,0);x < 0, y > 0\},\\
                 SwR            &=\{(x,y,0);x.y>0\},\\
                  EscR          &=\{(x,y,0);x> 0, y<0\}.
\end{array}
\label{regioes}\end{equation}

\subsection{Dynamic of $\varphi_Z$ when $Z$ is $\xi-$reversible} The next Lemma exhibit the
 subset of $\Omega_0^F(\delta)$ that is
$\xi$-reversible.

\begin{lemma}
Let $Z\in \Omega_0^F(\delta)$. If $XYh(0)=-YXh(0),$
$X^2h(0)=-Y^2h(0)$ then $Z$ is $\xi$-reversible.
\label{lema-reversivel}\end{lemma}

\dem. Straight forward computations.

$\fp$

The proof of Theorem $\ref{existenciadeciclos}$ item $(d)$ follows
by:

\begin{lemma}
Let $Z=(X,Y)\in \Omega_0^{\xi}(\delta)$. Then:
\begin{itemize}
    \item [$(i)$] 0 is a hyperbolic equilibrium saddle point of
    $\varphi_Z$ provided that $|XYh(0)|$ $>|X^2h(0)|$, see Figure \ref{tipo-primeiro-retorno};
    \item [$(ii)$] 0 is a hyperbolic equilibrium elliptical point of
    $\varphi_Z$ provided that $|XYh(0)|<|X^2h(0)|$.
\end{itemize}
\label{proposicao-reversivel}\end{lemma}


\dem. Consider the local normal form of $Z=(X,Y)$ given in the
Proposition $\ref{normal-forms}$. By the Lemma \ref{lema-reversivel}
we have: $XYh(0)=C_{XY},YXh(0)=C_{YX},X^2h(0)=C_X$ and
$Y^2h(0)=C_Y$.

The regions in $\Sigma$ are given in $(\ref{regioes})$. Let
$p_0=(x_0,y_0,z_0)$. The flows of the vector fields $X$ and $Y$ are:

\begin{equation}
\begin{array}{ll} \phi_X^{t}(p_0)       &=\left(x_0 + C_Xt,y_0 + C_{XY} t, z_0 + x_0t +
\frac{1}{2}C_X t^2\right)+O(t^2,t^2,t^3)\\\\
               \phi_Y^{t}(p_0)         &=\left(x_0 + C_{YX} t, y_0 + C_Y t, z_0+ y_0t+
               \frac{1}{2}C_Y t^2\right)+O(t^2,t^2,t^3).
\end{array}
\label{fluxo-dobra-dobra}\end{equation} If $p_0=(x_0,y_0,0)\in
\Sigma$, taking

\[
t_1(p_0)=-\frac{2 x_0}{C_X} \mbox{\,\,\, and \,\,\,}
t_2(p_0)=\frac{-4C_{XY} x_0+2C_X y_0}{C_X^2} \] we obtain
$\phi_X^{t_1}(p_0)$ $=(x_1,y_1,0)\in \Sigma$ and
$\phi_Y^{t_2}(x_1,y_1,0)\in \Sigma$. We define the return region
(the shaded region in the Figure \ref{tipo-primeiro-retorno}) by:

\[
R_0=[\{(x,y,0)\in \Sigma; x>0,y < 2C_{XY}  C_X^{-1} x\}\cap SwR].
\]

The return region $R_0$ represent the points $p=(x,y,0)\in\Sigma$
such that $t_i(p),i=1,2$ are positive. The first return map is
expressed by:

\begin{figure}[ht]
\epsfysize=2,3cm
\psfrag{On}{$\omega_-$}\psfrag{Op}{$\omega_+$}\psfrag{Vp}{$v_+$}\psfrag{R}{$R_0$}
\psfrag{Vn}{$v_-$}
\psfrag{X}{$x$}\psfrag{C}{$C$}\psfrag{Y}{$y$}\psfrag{1}{Elliptical
type}\psfrag{2}{Saddle type}\psfrag{M}{$\Sigma$}
\centerline{\epsfbox{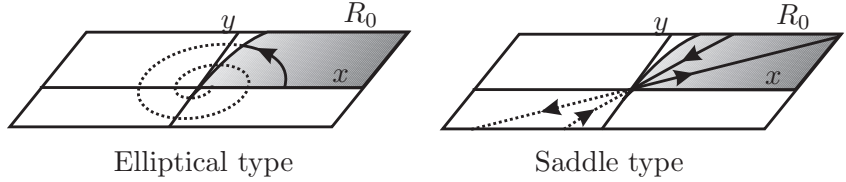}} \caption{Dynamics
of $\varphi_Z$.}\label{tipo-primeiro-retorno}
\end{figure}

\begin{equation}
\varphi_Z(p)=\phi_{Y}^{t_2}\circ\phi_{X}^{t_1}(p)=\left(\left(\frac{4C_{XY}^2-C_X^2}{C_X^2}\right)x
-\frac{2C_{XY}}{C_X}y,\frac{2C_{XY}}{C_X}x-y\right).
\label{primeiro-retorno-reversivel}\end{equation} The eigenvalues of
$\varphi_Z$ are $\lambda_{\pm}=-1+\frac{2(C_{XY}^2 \pm
\sqrt{C_{XY}^2(C_{XY}-C_X)(C_{XY}+C_X)})}{C_X^2}$. Besides, the
origin is a critical point of the saddle (respectively elliptical)
type of $\varphi_Z$ if $|C_{XY}|>|C_X|$(respectively
$|C_{XY}|<|C_X|)$.

$\fp$

\subsection{Existence of $k-$periodic orbits for $Z$} We investigate
the existence of periodic orbits of $Z\in \Omega_0^{\xi}(\delta)$ of
period $k$ arbitrary. We prove the items $(a),(b)$ and $(c)$ of
Theorem $\ref{existenciadeciclos}$:

\dem. By $(\ref{primeiro-retorno-reversivel})$, the equation
$\varphi_Z(x,y,0)=(x,y,0)$ is expressed by:

\[
\left\{\begin{array}{ll} \dis\frac{-C_X^2 + 4C_{XY}^2}{C_X^2}x
-\dis\frac{2C_{XY}}{C_X}y                      &=x\\\\
\dis\frac{2C_{XY}}{C_X}x -y                 &=y.
\end{array}\right.
\] The solutions are $C_{XY}=\pm C_X$. So, the periodic orbits pass
through the pair of straight lines $\{(x,y); x=\pm y\}$. Observe
that the straight line

\[
\mathcal{C}_{-1}=\{(x,y,0); x=-y\}\] is contained in $SlR\cup EscR$.
So we get a one-parameter family of pseudo periodic orbits.

On the straight line

\[
\mathcal{C}_1=\{(x,y,0);x=y\}=Fix\,(\xi)
\]we obtain a one-parameter family of periodic orbits. The expression of $\varphi_Z^n$ is given by:

\[
\begin{array}{cl}  \psi^n(x,\alpha x)=\varphi_Z^n(x,\alpha x)     &=
(\left[\frac{-C_X^2+4C_{XY}^2}{C_X^2}
-\dis\frac{2\alpha C_{XY}}{C_X}\right]^nx,\\\\
                                                     &\left[\frac{-C_X^2+4C_{XY}^2}{C_X^2}
                                                     -\frac{2\alpha C_{XY}}{C_X}\right]^{n-1}
                                                     \left[\frac{2C_{XY}}{C_X} -\alpha\right]x ).
\end{array}
\]

Solving $\psi^n(x,\alpha x)=(x,\alpha x)$, we obtain the system

\[
\left\{\begin{array}{lr}
\left[\dis\frac{-C_X^2+4C_{XY}^2}{C_X^2}-\dis\frac{2\alpha
C_{XY}}{C_X }\right]^n
                    &=1\\\\
                           \left[\dis\frac{-C_X^2+4C_{XY}^2}{C_X^2}-\dis\frac{2 \alpha C_{XY} }{ C_X }
                           \right]^{n-1} \left[\dis\frac{2C_{XY}}{C_X}-\alpha\right]
                           &=\alpha.
\end{array}\right.\]

We obtain directly that $C_X=\frac{2 \alpha C_{XY}}{1+\alpha^2}$.
Replacing this value in the first equation, we obtain
$\alpha^{-2n}=1$, ie, $\alpha=\pm 1$. So the conclusion of the
Theorem $\ref{existenciadeciclos}$ is straight forward.

$\fp$

\section{Pseudo periodic orbits}\label{pseudo-periodica}

\subsection{Existence of $k$-pseudo periodic orbits for $Z$} 
\begin{figure}[ht]
\epsfysize=2.7cm
\psfrag{X}{$X$}\psfrag{1}{$(a)$}\psfrag{2}{$(b)$}\psfrag{Y}{$Y$}
\psfrag{Sx}{$S_X$}\psfrag{Sy}{$S_Y$}
\psfrag{Rcu}{$SwR^+$}\psfrag{Rcd}{$SwR^-$}
\psfrag{F1}{$\varphi_1^{}$}\psfrag{Rd}{$SlR$}\psfrag{F2}{$\varphi_2^{}$}
\psfrag{P}{$p_0$}\psfrag{Re}{$EscR$}
\psfrag{R}{$r$}\centerline{\epsfbox{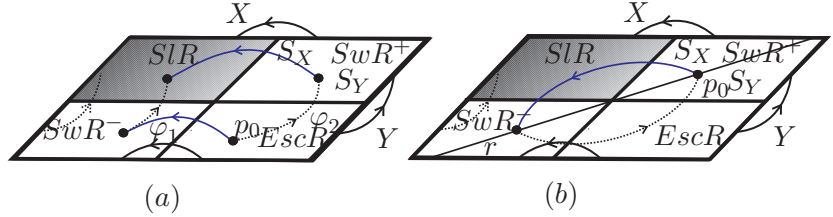}}
\caption{In $(a)$ is represented a pseudo periodic orbit. In $(b)$
is represented one of the periodic orbits of $Z\in
SL(\omegar)$.}\label{prop-orbitapseudoperiodica}
\end{figure}

A necessary
condition for the existence of pseudo periodic orbits is

\[
\varphi_1(x,y)=\varphi_2(x,y)
\] where
$\varphi_1(x,y)=(\phi_{Y}^{t_2}\circ \phi_{X}^{t_1})$ and
$\varphi_2(x,y)=(\phi_{X}^{t_4}\circ \phi_{Y}^{t_3})(x,y,0)$, with
$t_i\geq 0$, for $i=1, \dots,4$. See Figure
$\ref{prop-orbitapseudoperiodica}$.

The regions on $\Sigma$ are given in $(\ref{regioes})$. We define
the subset of $EscR\subset \Sigma$:

\[
RS=\left\{(x,y,0)\in EscR; \dis\frac{2Y(Xh)(0)}{Y^2h(0)}y < x <
\dis\frac{X^2h(0)}{2X(Yh)(0)}y\right\}.
\]

\begin{lemma}
Let $Z\in SL(\Omega_0^F(\delta))$. If
$X^2h(0)<0,Y^2h(0)>0,XYh(0)=-YXh(0),X^2h(0)=-Y^2h(0)$ and
$XYh(0)YXh(0)=X^2h(0)Y^2h(0)$ then $RS=\emptyset$.
\label{lema-regiao}\end{lemma}

\dem. In fact,

\[
\dis\frac{2YXh(0)}{Y^2h(0)}y < \dis\frac{X^2h(0)}{2XYh(0)}y
\Longleftrightarrow 4X^2h(0)Y^2h(0)>X^2h(0)Y^2h(0).
\] But this is a contradiction since
$X^2h(0)<0$ and $Y^2h(0)>0$, by hypothesis.

$\fp$

We recall the notations:

\[
\mathcal{C}_{\alpha}=\{(x,\alpha x,0)\in \Sigma; \mbox{ with }\alpha
\in \R\}.
\]

Now we conclude the proof of Theorem
$\ref{proposicao-orbita-pseudo-periodica}$:

\vs

\textbf{Proof of Theorem
$\ref{proposicao-orbita-pseudo-periodica}$}: Consider the regions in
$\Sigma$ given by $(\ref{regioes})$. The flow given in
$(\ref{fluxo-dobra-dobra})$. So, the expressions of $\varphi_1(x,y)$
and $\varphi_2(x,y)$ are:

\begin{equation}
\begin{array}{ll}   \varphi_1(x,y)             &=\left(\left(-1+\dis\frac{4C_{XY} C_{YX}}{C_XC_Y}\right)x
-\dis\frac{2C_{YX}}{C_Y}y,-y+\dis\frac{2C_{XY}x}{C_X} \right),\\\\
                    \varphi_2(x,y)              &=\left(-x+\dis\frac{2C_{YX}}{C_Y}y,-\dis\frac{2C_{XY}}{C_X}x
                    +\left(\dis\frac{4C_{XY}C_{YX}}{C_XC_Y}-1\right)y\right),
\end{array}
\label{aplicacoes-primeiro-retorno}\end{equation}where
$t_1(x,y)=-\frac{2x}{C_X},t_2(x,y)=-\frac{2}{C_Y}(-\frac{2C_{XY}}{C_X}x+y),
t_3(x,y)=-\frac{2}{C_Y}y$ and $t_4(x,y)$
$=-\frac{2}{C_X}(x-\frac{2C_{YX}}{C_Y}y)$. As $X^2h(0)<0$ and
$Y^2h(0)>0$ we get $t_1>0$ and $t_3>0$, respectively.

To prove item $(a)$ we need to solve $\varphi_1(p)=\varphi_2(p)$
with $p=(x,y,0)\in \Sigma$. Initially we suppose that $p\in EscR$.
From $(\ref{aplicacoes-primeiro-retorno})$ and solving
$\varphi_1(x,y)=\varphi_2(x,y)$ we obtain:

\begin{itemize}
\item [$(i)$] $(X(Yh))(0)(Y(Xh))(0)$ $=(X^2h)(0)(Y^2h)(0)$;
\item [$(ii)$] $(x,y)\in \mathcal{C}_{\lambda^*}$, where $\lambda^*=\frac{X^2h(0)}{XYh(0)}$.
\end{itemize}

As $\mathcal{C}_{\lambda^*}\subset EscR$ we obtain that $XYh(0)>0$.
Replacing this condition in $(i)$ we get $YXh(0)<0$.

In this way, we consider the diffeomorphisms $\varphi_1$ and
$\varphi_2$ restricted to $RS\subset \Sigma$. By Lemma
$\ref{lema-regiao}$, $RS=\emptyset$.

Therefore, for $Z\in SL(\Omega_0^F(\delta))$ there does not exist
pseudo periodic passing by $p\in EscR$.

If $p\in SwR$ follows by item $(c)$ of Theorem
$\ref{existenciadeciclos}$ the result. Observe that in this case,
some values of the time $t_i,i=1,\dots,4$ are equal to 0.

Let us prove item $(b)$. For $Z\in\Omega_0^P(\delta)$, the
conditions $(i),(ii)$ and $\varphi_1=\varphi_2$ are satisfied since
$\mathcal{C}_{\lambda^*}\subset SwR$. The dynamics of $Z\in
\Omega_0^P(\delta)$ is illustrated in Figure
$\ref{prop-orbitapseudoperiodica}$.

\textbf{Claim}: ``If $p\in \mathcal{C}_{\lambda^*}\subset SwR^+$ and
$Z\in \Omega_0^P(\delta)$ (respectively $p\in
\mathcal{C}_{\lambda^*}\subset SwR^-$ and $Z\in \Omega_0^P(\delta)$)
then $\varphi_1(p)=p$ (respectively $\varphi_2(p)=p$). In other
words, on $\mathcal{C}_{\lambda^*}$ there is a family of periodic
orbits''.

In fact, observe that if $Z\in \Omega_0^P(\delta),p\in
\mathcal{C}_{\lambda^*}\subset SwR$ then $Z$ satisfies the
conditions $(i),(ii)$ and $t_i\geq 0$, for $i=1\dots, 4$. That is,
$\varphi_{1}(p)=p$ $(p\in \mathcal{C}_{\lambda^*} \cap SwR^+)$ or
$\varphi_{2}(p)=p$ $(p\in \mathcal{C}_{\lambda^*} \cap SwR^-)$.
Besides, for $Z\in \Omega_0^P(\delta)$ we have
$\mathcal{C}_{\lambda^*}\subset (SwR^+ \cup SwR^-)$. If $p\in
[\mathcal{C}_{\lambda^*}\cap SwR^-]$ then
$\varphi_1(p)=\phi_{Y}^{t_2}\circ \phi_{X}^{t_1}(p)=-p$ and
$t_1(x,y)=0$.

$\fp$

\medskip

\noindent {\textbf{Acknowledgments.}} The first and third authors
wish to thank Brazil-France cooperation agreement. The second author
is partially supported by FAPESP--BRAZIL grant number 07/56163-4.
The third author is partially supported by FAPESP--BRAZIL and
CNPq--BRAZIL. This work is partially realized at UFG/Brazil as a
part of project numbers 35796 and 35797.

\end{document}